\documentclass{amsart}
\usepackage{amssymb,amscd}

\newtheorem{definition}{Definition}
\newtheorem{theorem}{Theorem}
\newtheorem{corollary}{Corollary}

\begin{document}

\title{RSA cryptosystem for rings with commuting ideals} 
\author{Nasrutdinov M.F., Tronin S.N.} 


\title{RSA cryptosystem for rings with commuting ideals}

%
%
%

\begin{abstract} 
	This article presents a generalization of the RSA cryptosystem for rings with commuting ideals. An analogue of the Euler function for ideals and the concept of an RSA-ideal are defined. An analog of a cryptosystem for the ring with commuting ideals is formulated and a description of the RSA-ideals for which this is possible is obtained.
\end{abstract}


\keywords{Algebraic cryptography, rings with commuting ideals, RSA cryptosystem, Euler's
	function for ideals.} 

\maketitle

\section{Introduction} \label{ }

In \cite{TP} the authors generalized the well-known RSA cryptographic algorithm to Dedekind rings.  They replaced natural numbers with ideals of commutative Dedekind rings. After some time it became clear that a similar idea could be realized for rings that are not necessarily commutative or Dedekind ring. This paper describes the initial stage of the implementation of this idea. 

RSA is a public-key cryptosystem that is widely used for secure data transmission. Ron Rivest, Adi Shamir and Leonard Adleman described the algorithm in 1977 \cite{RSA}.  This cryptosystem is the object of intense study up to the present time.

Recall the idea of the original RSA \cite{RSA}. The receiver of a message publishes the number $n$, where $n=pq$ is the product of two different secret large primes and the encryption exponent $e$ satisfying the condition that the greatest common divisor $gcd (e, \varphi(n))=1$. Here $\varphi(n)$ is the Euler's function. 

Let $m$ be a bit string which  we want to encrypt. Let's represent $m$ as a natural number from $0$ to $n$. The sender of the letter computes the ciphertext represented by the number (or bit string) $c=m^e \pmod {n},$ which is the encrypted message. The receiver of the message selects $1<d < \varphi (n)$ from the condition $ed=1+ \varphi(n)t.$ Then the receiver recovers the original message by raising the ciphertext to the power of $d$.  Because $m^{ed} =m \pmod{n}$.

The secret key is a triple $p$, $q$, $\varphi (n)=(p-1)(q-1)$. Without knowing $p$ and $q$ it is very difficult to calculate $\varphi (n)$, without $\varphi (n)$ it is very difficult to calculate $d$, i.e. to decrypt the message. The security of the RSA cryptosystem is based on the factorization complexity $n$. Details can be found in the \cite{Song} book.

The security of the RSA cryptosystem is based on two mathematical problems: the problem of factoring large numbers and the RSA problem.   But the rise of quantum computing  becomes the threat to modern cryptography. Shor's quantum algorithm \cite{Shor}, in particular, provides a large theoretical speedup to the brute-forcing capabilities of attackers targeting many public-key cryptosystems such as RSA (see also \cite{Song}).  Therefore, there is a practical need to find new cryptographic algorithms, attacks on which will be difficult even for quantum computers. We investigate the possibility of using commutative rings that are not principal ideal rings or non-commutative rings to generalize the RSA cryptosystem. We show that a meaningful generalization of RSA is possible for a fairly wide class of associative rings with identity.

Let us briefly describe the content of the paper. In section 2 we consider the special  class of rings, namely rings witch commuting ideals (CI-rings). These rings are quite convenient for generalizing the RSA cryptosystem.  In section 3 we generalize Euler's function for ideals of CI-rings.  In section 4 we construct an analog of the RSA algorithm for CI-rings and we define the concept of an RSA-ideal. The main result of this paper is  
a description of the RSA-ideals.  In fact, to generalize the RSA cryptosystem we should  use only intersection of maximal ideals.  

The results of this research were announced in \cite{TrNas}.

\section{Rings with commuting ideals}

All rings are assumed to be associative and with nonzero identity element. We will consider rings with commuting ideals. As will be shown below, this class of rings allows generalization of the RSA cryptoscheme and contains both commutative and non-commutative rings.


\begin{definition}
A ring $R$ is said to be a CI-ring (ring with commuting ideals) if  $AB = BA$ for all ideals $A, B$ of $R$. 
\end{definition}

As usual  the product $AB$ of two ideals $A$ and $B$ is the set of all finite sums of elements of the form $ab$ with $a \in A$ and $b \in B$. 

CI-rings have the following easily proved properties: 

\begin{theorem}{\rm \cite{Arm}}
	\begin{enumerate}
		\item A homomorphic image of a CI-ring is a CI-ring.
		\item A finite direct product of CI-rings with unity is a CI-ring.
		\item If $R$ is a CI-ring with unity, then for each positive integer $n$, the ring
		of $n  \times n$ matrices over $R$ is also a CI-ring.		
		\item  Let $R$ be a CI-ring and $e=e^2 \in R$ is the idempotent. Then $eRe$ is a CI-ring. 
	\end{enumerate}	
\end{theorem}

{\bf Proof}
The first two statements are obvious. 

To prove the third statement, note that every ideal of $M_n(R)$ is of the form $M_n(A)$ where $A$ is an ideal of $R$. 

To prove the last statement we use the fact that $RR=R$. Suppose $A,B$  are ideals of $eRe$. We have  $$AB = eRe A eRe eRe B eRe = e(Re A eR) (R e eRe B eR) e = e (R e eRe B eR)  (Re A eR) e = B A.$$

CI-rings were introduced by Armendariz and Heatherly in short proceedings thesis \cite{Arm}. At the moment, we do not know the literature where the theory of such rings would develop. There are some examples of such rings in \cite{Jacobson_en} and \cite{Cohn}. In \cite[Chapter 4]{Tug} modules over CI-rings are considered.

\vskip 0.3cm 

{\bf Examples of CI-rings} 

\vskip 0.3cm 

1) Commutative rings. 

2) Principal ideal domains. Recall that $R$ is a  principal ideal domain if $R$ is a domain of integrity and every  one-sided ideal (right or left) is a principal one-sided ideal (right or left respectively).  These rings are considered in the third chapter of \cite{Jacobson_en}. 

Obviously, the rings of integers is a principal ideal domain. There are other examples of this ring: 

2.1) The subring of Hamilton's quaternion algebra consisting of quaternions $\alpha_0+i \alpha_1+ j \alpha_2+ k \alpha_3$,
 where the $a_i$ are either all rational integers or all halves of odd integers. This subring is called the ring of Hurwitz quaternions. 
 
2.2) Let $K$ be any skew field with an endomorphism $\alpha$ and an $\alpha$-derivation $\delta$. Then the skew polynomial ring $K[x; \alpha, \delta]$ is a principal ideal domain whenever $\alpha$ is an automorphism \cite[Theorem 1.3.2]{Cohn}.

3) Leavitt path algebras of directed graphs \cite[Theorem 4.2]{Rang}.    

For arbitrary rings the famous Chinese remainder theorem holds. 

\begin{theorem}{\rm \cite[Theorem 18.30]{Faith}}
If $A_1,  \dots, A_n$ are finitely many ideals of a ring $R$, then the following conditions are equivalent: 
\begin{enumerate}
	\item The canonical map $h: R/\bigcap\limits_{i=1}^n A_i \to \prod\limits_{i=1}^n R/A_i$  with $h(r) = (r+A_1, r+A_2, \ldots, r+A_n)$ is an isomorphism.
	\item For any set $x_1, x_2, \ldots, x_n$  of elements of $R$ the system of congruences   $X = x_i (\pmod A_i)$ has a solution $x \in R$. 
	\item The ideals $A_1,  \dots, A_n$ are comaximal in pairs, that is, $A_i +A_j =R$ whenever $i \not = j$.
\end{enumerate}	
\end{theorem} 

For CI-rings we can prove additionally.

\begin{corollary}
	Let $R$ is a CI-ring and $A_1,  \dots, A_n$ are comaximal in pairs. Then $$ \bigcap\limits_{i=1}^n A_i = \prod\limits_{i=1}^n A_i.$$
\end{corollary}

{\bf Proof} We use induction on $n$. 
	
Let $n = 2$. Suppose $A,B$  are comaximal ideal and $AB=BA$. We shall prove that $A \bigcap B = AB$. It is sufficient to show that $A \cap B \subset AB$. From $R=A+B$ we have   $1=a+b$ for some $a \in A$ and $b \in B$. Then for all $x \in A \cap B$ we have 
$x=x(a+b) \in AB +BA =AB$.

Now, suppose that $n > 2$ and that the result has been proved for smaller values of
$n$. For pairwise comaximal ideals $A_1,  \dots, A_n, A_{n+1}$ of $R$ we denote by $B  = \bigcap\limits_{i=1}^n A_i = \prod\limits_{i=1}^n A_i$. We claim that $$\prod\limits_{i=1}^{n+1} A_i= BA_{n+1}= 
B\bigcap A_{n+1}=\bigcap\limits_{i=1}^{n+1} A_i.$$

By the Chinese remainder theorem there is $x \in R$ such that  $x = 0 \pmod{A_i}$  for all 
$i=1,\ldots, n$ and $x=1 \pmod {A_{n+1}}$.  
Then $x \in \bigcap\limits_{i=1}^n A_i = B$ and $1-x \in A_{n+1}$. Thus, for any $r \in R$ we have $$r=rx+ r(1-x) \in B+ A_{n+1}.$$ Therefore $B, A_{n+1}$ are comaximal ideals and $BA_{n+1}=B\cap A_{n+1}$.

{\bf Note.} Corollary 1 will allow us to decompose RSA ideals of CI-rings into product of maximal ideals and to introduce the Euler function for ideals with suitable properties.

\section{Euler's function for ideals}

In \cite{TP} the Euler's function analog was defined for ideals of Dedekind rings. Note that RSA cryptosystems with Dedekind rings were also considered in \cite{Kondratyonok}. We  define the Euler's function in the same way.   

\begin{definition}
Let $A$ be the ideal of a ring $R$ and $|U(R/A)| < \infty$ where $U(R/A)$ is a group of units of $R/A$.  
The function $$ 
\varphi(A)=|U(R/A)|.
$$ is called the Euler's function of the ideal $A$. 
\end{definition}

{\bf Example 1:} If $R= \mathbb Z$ then $A=(n)= \mathbb Zn$ and $\varphi (A)= \varphi (n)$ in the usual sense.

{\bf Example 2:} If $R$ is the ring of Gaussian integers . Then 
$$
\varphi (m)=|m|^2 \mathop{\prod}\limits_{p|m}\,\left(1-\frac{\displaystyle 1} {\displaystyle |p|^2}\right).
$$

{\bf Example 3:} If  $R=\mathbb{Z}[\sqrt{k}]$ is a ring of the quadratic integers. Then 
$$
\varphi (m)=|\nu(m)| \mathop{\prod}\limits_{p|m}\,\left(1-\frac{\displaystyle 1} {\displaystyle N(m)}\right),
$$ 
where $N(a+b\sqrt{k})=a^2-kb^2$ is norm of element of $m=a+b\sqrt{k}\in R$ and  $\nu(m)=|R/(m)|$.

{\bf Example 4:} If $R$ is a polynomial ring over Galua field $GF(q)$ then
$$\varphi (m)=q^{deg (m)} \mathop{\prod}\limits_{p|m}\,(1-q^{-deg (p)}).$$

\begin{theorem} Let $R$ be a CI-ring.  For finite set of pairwise comaximal ideals $A_1,  \dots, A_n$ of the ring $R$ we have $$\varphi (A_1 A_2 \ldots A_n)= \varphi (A_1) \varphi (A_2) \ldots  \varphi (A_1).$$
\end{theorem}

{\bf Proof}
By the Chinese remainder theorem and Corollary 1 we have $$R/ \prod\limits_{i=1}^{n} A_i = 	\prod\limits_{i=1}^n R/A_i.$$ 
	$$ \varphi_R(\prod\limits_{i=1}^{n} A_i)= |U (R/ \prod\limits_{i=1}^{n} A_i)| = 
	U(\prod\limits_{i=1}^n R/A_i) = 
	\prod\limits_{i=1}^{n} |U (R/A_i)| = \varphi (A_1) \varphi (A_2) \ldots  \varphi (A_1).$$

\section{RSA Cryptosystem for CI-rings and RSA-ideals}

The RSA algorithm uses natural numbers.The security of RSA relies on the practical difficulty of factoring the product of two large prime numbers. We replace the ring $\mathbb Z$ in the scheme by some CI-ring $R$ and the integers by  ideals of the ring $R$. Additionally, we assume that the Euler's function for the ideals must exist and be difficult to compute. 

Suppose $2< |R/A|<\infty.$ If $\gcd(e,\varphi(A))=1$ and $ed=1+ \varphi(A)t$ for some natural numbers $e,d$ then we need $m^{ed}\equiv m \pmod{A}$ holds for all $m \in R.$ The encryption and decryption algorithm in case of the RSA cryptosystem for CI-rings is similar to the classical case. Let's describe it. 

Suppose Bob wants to send messages to Alice that only she can decrypt. Alice chooses two ideals $M_1 \neq M_2\subset R$. Then she calculates the ideal $A =M_1 \cdot M_2$. Also Alice chooses an encryption exponent $e$ that satisfies the condition $\gcd(e,\varphi(A))=1.$ So, Alice's public key is the pair $(A, e)$ and she can keep it in the public domain. After that, Alice calculates the decryption exponent $d$ from
conditions $ed=1+\varphi({A})t$. Alice's secret key is the triple $(d, M_1, M_2)$, which she keeps secret.

Suppose Bob wants to encrypt a message for Alice. He represents the message as an element $m \in W$. The ciphertext $c$ is obtained by raising the message to a power equal to the open encryption exponent and taking the remainder modulo $A$
$$
c=m^e \pmod A.
$$

Alice can decrypt the ciphertest $c$ and get the original message as follows:
$$
 m=c^d \pmod{A}. 
$$

The message $ m=c^d$ is represented by an element from the set $W$. $W$ is a complete system of remainders modulo $A$. Thus $|(r+A) \cap W|=1$ for each $r \in R.$ There is a bijection $\phi : W \leftrightarrow R/A$. At the same time elements from $W$ should be convenient for encoding into a bit string.

In order for the cryptoscheme would work we need the condition $m^{ed} = m \pmod{A}$ holds. We introduce the following definition.

\begin{definition} Suppose $A$ be a ideal of CI-ring $R$, $R/A$ is a finite ring,  $\varphi(A)>2$, $e, d$ are natural numbers such that $1<e,d< \varphi (A)$ and  $\gcd(e,\varphi (A))=1,$ $ed\equiv 1+\varphi (A)t$. 
	
We shall say that ideal $A$ is said the RSA-ideal if $x^{ed} \equiv x \pmod {A}$ for all   $x \in R$.
\end{definition}

The following theorem is a direct consequence of Jacobson's theorem  \cite[Theorem 3.1.2]{Herst} on the commutativity of rings with the property $x^{n(x)} =x$.

\begin{theorem} 
	Let $A$ be an ideal of a CI-ring $R$ then $R/A$  is a commutative ring.
\end{theorem}

Our main result is the following.

\begin{theorem} Let $A$ be an ideal of a CI-ring $R$.  The ideal $A$ is 
RSA-ideal if and only if $A=M_1 M_2 \ldots M_k$ where $M_i$ is maximal ideal for any $i=1,\ldots,k$. 
\end{theorem}

{\bf Proof}

$(\Rightarrow)$  Consider the factor ring $R/A$. Denote by $s=ed = 1 \pmod{\varphi(A)}$. Then for any $x \in R/A$ we have $x^s =x$. 


Let us show that the quotient ring $R/A$ has no nilpotent elements. Assuming the converse, let $x \not = 0$ and $n$ be the smallest number for which $x^n=0$. 

If $n < s$ then $x=x^s=x^{s-n}x^n =0$ and this contradicts to the choice of $x$. 

If $n>s$ then we divide $n$ with remainder $n=qs+r$, $q>0$, $0 \leq r<s$. Then $0 = x^{n} = x^ {qs+r} = x^{sq}x^r=x^{q+r}$ and we have a contradiction with the minimality of $n$.


By the Wedderburn-Artin theorem, $R/A$ is isomorphic to a direct sum of fields. Indeed, $R/A$ is a finite commutative ring without nilpotent elements. Hence its Jacobson radical is equal to zero \cite[Theorem 1.3.1]{Herst} and $R/A$ is the Artinian semisimple ring. Thus, $R/A$  is isomorphic to the direct sum of matrix rings over division rings \cite[Theorem 1.4.4 and Theorem 2.1.6]{Herst}.

Since the matrix ring of order $\geq 2$ is not commutative it follows that $R/A \cong F_1 \oplus F_2 \oplus F_k$ where $F_i$ is a field for each $i$.  

Consider an epimorphism $\theta : R \to R/A$ and projections $\pi_j: R/A \to F_j$. Then $M_j = \ker \pi_j \theta$ are maximal ideals. By the Chinese remainder theorem $R/ \bigcap\limits_{j=1}^k M_j = \prod\limits_{i=1}^n R/A_i$. Therefore, $$\bigcap\limits_{j=1}^k \ker \pi_j \theta = \bigcap\limits_{j=1}^k M_j = M_1 M_2 \ldots M_k.$$

We claim that $A=\bigcap\limits_{j=1}^k \ker \pi_j \theta$.  Clearly, $A = \ker \theta  \subset  \bigcap\limits_{j=1}^k \ker \pi_j \theta$. Conversely, if $r \in R \not \in A$ and $r \in \bigcap\limits_{j=1}^k \ker \pi_j \theta$ then $\theta(r) = r+A  \not =0$ but 
$(\pi_1 \theta(r), \pi_2 \theta(r), \ldots, \pi_k \theta(r)) = (0,0, \ldots, 0)$. It is impossible because  $R/A$ and $\prod\limits_{j=1}^n R/M_j$ are isomorphic.

$(\Leftarrow)$ Let $A = M_1 M_2 \ldots M_k$ be the product of maximal ideals,
$s=ed=1+t \cdot \varphi(A)=1+ t\cdot \varphi(M_1) \varphi(M_2) \ldots \varphi(M_k)$. Then $F_j=R/M_j$ are finite skew fields (fields, due to the commutativity of finite skew fields) and by the Chinese remainder theorem $R/A \cong F_1 \oplus F_2 \oplus \ldots \oplus F_k$.

Then $U(R/A) = U(F_1^*) \times U(F_2^*) \times \ldots \times U(F_k^*)$ where 
$F^{*} = F \setminus \{ 0 \}$. If $x \in A$ then $x^s \in A$ and $x^s = x \pmod{A}$.

If $x \not \in A$ then the image of $x$ in $R/A$ can be represented as $$\overline{x} = (x_1, x_2, \ldots, x_k) \in F_1 \oplus F_2 \oplus \ldots \oplus F_k.$$ It suffices to show that $x_j^{s}=x_j$. 

If $x_j=0$ then the assertion is trivial. If $x_j \not =0$, then $x_j \in F_j^* = U(F_j)$ and by the Lagrange theorem on the element orders in the group $x_j^{\varphi(M_j)}= 1$. Recall that $\varphi(M_j) = |U(F_j)| $.
Thus, $$x_j^s = x_j^{1+ t\cdot \varphi(M_1) \varphi(M_2) \ldots \varphi(M_k) } = x_j.$$

\begin{corollary}
	In the definition of the RSA-ideal $A$ of ring $R$ any number $e$ coprime to $\varphi (A)$ is allowed. 
\end{corollary}

\section{Conclusions and open problems}

From a practical point of view, the ring for building an RSA cryptosystem should not be complicated. This means that the elements of this ring can be quite conveniently represented in the computer's memory, and operations with such elements are easily computable.

At the same time, to implement the RSA analogue protocol, it is necessary that the task of finding the secret key $d$ from public keys be a computationally difficult task.


It follows from the description of RSA ideals given in the theorem that the security  of a hypothetical encryption protocol will depend both on the complexity of calculating $\varphi (A)$ and on the complexity of the problem of decomposing this ideal into a product of maximal ideals. The uniqueness of such a decomposition is not assumed. The RSA ideal in this scheme will play the role of one of the public keys.

The results of our research  give a hint in which direction we can look for non-commutative rings suitable for application in practical cryptography. An open problem is if there are good rings apart from $\mathbb Z$.


The ring of Hurwitz quaternions if we replace of natural numbers by ideals does not provide any advantages compared to $\mathbb Z$. This follows from the description of the ideals of this ring \cite[Theorem~211]{Redei}. However, in this case, there is a fairly meaningful description of simple elements which are not directly related to the structure of ideals \cite[Chapter 5]{Conway}. So it is not worth completely excluding Hurwitzian quaternions from consideration yet.

The next candidates for research may be non-commutative principal ideal domains. They are considered fairly simple.  In \cite{Jacobson_en} it shows that the ideals of a non-commutative principal ideal ring commute. 

Another well-studied class of rings with commuting ideals are the skew polynomial rings \cite{Good}. Finally, it is of interest to find out which of the finite non-commutative rings may be suitable for use in the ring analogs of the RSA cryptosystem. To begin with, we can pose the question of characterizing finite non-commutative rings with commuting ideals.

Another interesting question from the point of view of ring theory is the following: Whether it is possible to weaken the conditions in the definition of RSA-ideals, for example, to weaken the condition that the quotient ring is finite.

\end{document}